\documentclass{amsart}

\usepackage{amssymb}

\usepackage[all]{xy}

\newtheorem{thm}{Theorem}

\newtheorem{cor}[thm]{Corollary}

\newtheorem{prop}[thm]{Proposition}

\theoremstyle{definition}
\newtheorem{defn}[thm]{Definition}

\newtheorem{say}[thm]{}
\newtheorem{exmp}[thm]{Example}

\newtheorem{rem}[thm]{Remark}          

\newtheorem*{ack}{Acknowledgments}      

\newtheorem{defn-thm}[thm]{Definition--Theorem}  
\newtheorem{defn-lem}[thm]{Definition--Lemma}  

\theoremstyle{remark}

\renewcommand{\c}[0]{{\mathbb C}}

\renewcommand{\a}[0]{{\mathbb A}}

\newcommand{\p}[0]{{\mathbb P}}

\newcommand{\map}[0]{\dasharrow}
\newcommand{\qtq}[1]{\quad\mbox{#1}\quad}

\newcommand{\im}[0]{\operatorname{im}}

\newcommand{\sing}[0]{\operatorname{Sing}}

\newcommand{\onto}[0]{\twoheadrightarrow}

\newcommand{\mor}[0]{\operatorname{Mor}}

\def\into{\DOTSB\lhook\joinrel\to}

\def\loccoh#1.#2.#3.#4.{H^{#1}_{#2}(#3,#4)}

\DeclareMathAlphabet{\mathchanc}{OT1}{pzc}%
                                {m}{it}

\newcommand{\chain}[0]{\operatorname{Chain}}

\begin{document}
\bibliographystyle{amsalpha}


\title[Lefschetz property]{The Lefschetz property for families of curves}
\author{J\'anos Koll\'ar}

\maketitle

By the Lefschetz hyperplane theorem, if $X\subset \p^N$ is a smooth,
projective variety and $C:=X\cap L$ is a
positive dimensional   intersection of $X$
with a linear subspace, then the natural map
$$
\pi_1(C)\to \pi_1(X)\qtq{is surjective.}
$$
The same conclusion holds if $X$ is quasi projective, but here
$C$ has to be an intersection of $X$
with a linear subspace in {\em general position.}
The aim of this note is to study families of curves $\{C_m:m\in M\}$
 that satisfy this Lefschetz--type property.
Such results were proved
 in the papers \cite{MR1786496, MR2011744}. 
Arithmetic applications are given in \cite{MR1786496, MR1745009, MR2019976}
and \cite{k-hom} studies this question
for homogeneous spaces. Related results are in 
\cite{MR2566998, MR2426347, MR2894633}.

\begin{defn}\label{leff.prop.defn}
A {\it family of schemes} over  a  normal variety $X$ over $\c$ is a diagram
$$
M\stackrel{p}{\leftarrow} C_M\stackrel{u}{\to} X.
\eqno{(\ref{leff.prop.defn}.1)}
$$
Our main interest is in families where $M$ is  irreducible
and $p$ is flat with irreducible  fibers.
The  family  (\ref{leff.prop.defn}.1)
satisfies the {\it Lefschetz property} if the following holds.
\begin{enumerate}\setcounter{enumi}{1}
\item[]   For every Zariski open  subset $\emptyset\neq X^0\subset X$
there is a Zariski open  subset $\emptyset\neq M^0\subset M$ such that, 
for every  $m\in M^0$, the induced map
$$
u(X^0,m)_*:\pi_1\bigl(C_m\cap u^{-1}(X^0)\bigr)\to \pi_1(X^0)\qtq{is surjective.}
$$
\end{enumerate}
We say that  (\ref{leff.prop.defn}.1)
satisfies the {\it weak Lefschetz property} if there is a constant
$N$  (independent of $X^0$) such that, for suitable choice of $M^0$,  
the image of  $u(X^0,m)_* $  
has index at most $N$ in  $\pi_1(X^0) $.

{\it Notes.} 
We ignore the base point since the surjectivity of the maps
between the fundamental groups of connected schemes 
does not depend on the choice of a base point.
The Lefschetz properties over arbitrary base fields are considered in
(\ref{leff.prop.defn.gen}).
\end{defn}

Our main Theorem \ref{main.thm} is somewhat technical, though I believe
it to be essentially optimal. 
The original arguments of  \cite{MR1786496, MR2011744}
need  high degree {\em very free} rational curves.
 By contrast, the current proof
frequently works for the lowest degree {\em free} curves.
As an illustration,
a simple yet  nontrivial example is given by lines on
hypersurfaces.

\begin{cor} \label{lines.on.hyp.cor}
Let $X\subset \p^n$ be a smooth hypersurface of degree $d$
over $\c$. The following are equivalent.
\begin{enumerate}
\item The family of lines has the Lefschetz property.
\item The family of lines has the weak Lefschetz property.
\item  $d\leq n-2$  or $X$ is a line in $\p^2$.
\end{enumerate}
\end{cor}

Let us start with some situations when the
Lefschetz property  fails.

\begin{exmp} \label{3.key.exmps}
Let  $M\leftarrow C_M\stackrel{u}{\to} X$ 
be  a flat, irreducible family of
irreducible varieties.

(\ref{3.key.exmps}.1) Assume that  $u$ is not dominant. Then 
any  $X^0\subset X\setminus \overline{u(C_M)}$
with infinite fundamental group shows that the weak
Lefschetz property does not hold. 

(\ref{3.key.exmps}.2) Assume that there is an open subset $X^*$ and a 
dominant morphism to a 
positive dimensional variety $q:X^*\to Z^*$
such that every $X^*\cap C_m$ is contained in a fiber of $q$ for general
$m\in M$. 
(We say that $X$ is {\it  generically $C_M$-connected}
if there is no such map  $q:X^*\to Z^*$;
see (\ref{chains.say.8}) for a better definition.)

Let $Z^0\subset Z^*$ be an open subset 
and $X^0:=q^{-1}(Z^0)$. Then
$$
\im\bigl[\pi_1\bigl(C_m\cap u^{-1}(X^0)\bigr)\to \pi_1(X^0)\bigr]
\subset \ker\bigl[ \pi_1(X^0)\to \pi_1(Z^0)\bigr].
$$
If $Z^0$ has   infinite fundamental group
then the  weak
Lefschetz property fails.

(\ref{3.key.exmps}.3) Assume that 
 $u:C_M\to X$ does not have geometrically irreducible generic fiber.
Then there is  a nontrivial Stein factorization 
$u:C_m \stackrel{w}{\map} Y \stackrel{v}{\to} X$
where $v$ is finite, generically \'etale of degree $>1$ 
and $w$ has geometrically irreducible generic fiber
(see \cite[Lem.9]{MR2011744} for the non-proper variant used here).
Let $X^0\subset X$ be an open set such that $v$ is \'etale
over $X^0$ and $Y^0:=v^{-1}(X^0)$. For general $m\in M$,
the induced map  $C_m\to X$ factors through $Y$, hence
$$
\im\bigl[\pi_1\bigl(C_m\cap u^{-1}(X^0)\bigr)\to \pi_1(X^0)\bigr]
\subset \im\bigl[ \pi_1(Y^0)\to \pi_1(X^0)\bigr]
\subsetneq\pi_1(X^0).
$$
In this case the Lefschetz property fails but the weak variant could hold
with $N=$ the number of geometric irreducible components
of the generic fiber of $u$.
More generally, we see that the  weak 
 Lefschetz property for $M\leftarrow C_M\stackrel{u}{\to} X$ 
is equivalent to the  weak 
 Lefschetz property for $M\leftarrow C_M\stackrel{w}{\to} Y$.
 The advantage is that $w:C_M\to Y$ has geometrically irreducible generic fiber.

(\ref{3.key.exmps}.4) An extreme case of the above is when
$u:C_M\to X$ is generically finite. Then $w:C_M\to Y$ is birational
and (\ref{3.key.exmps}.2) applies to $p:C_M\to M$. Thus
the   weak 
 Lefschetz property does not hold for $M\leftarrow C_M\stackrel{u}{\to} X$.
(A trivial exception is when $M$ is a single point,
giving the  case $X=\p^1$  in (\ref{lines.on.hyp.cor}.3).)

(\ref{3.key.exmps}.5) A difficulty in using the reduction method
of (\ref{3.key.exmps}.3) is that being generically $C_M$-connected
changes as we pass from $X$ to $Y$. A typical example is the following.
 
For some $n\geq 2$ set $X=S^2\p^n\setminus (\mbox{diagonal})$.
 Its universal cover is 
$\tilde X=\p^n\times \p^n\setminus (\mbox{diagonal})$.
Let $\tilde u:\tilde C_L\to \tilde X$ be the family of lines 
that are  contained in 
some $\p^n\times \{\mbox{point}\}$ and $u:C_L\to X$ the 
corresponding family of lines in $X$. Note that
$X$ is generically $C_L$-connected but $\tilde X$ is not
generically $\tilde C_L$-connected.

Since each point in $X$ has 2 preimages in $\tilde X$,
each fiber of $u$ has 2 irreducible components.

For an open set $W\subset \p^n$ let
$X_W\subset X$ denote the image of $W\times W$. Then
there is an extension
$$
1\to \pi_1(W)+\pi_1(W)\to \pi_1\bigl(X_W\bigr)\to \{\pm 1\}\to 1
$$
and for any line $C_m$, the image of
$\pi_1\bigl(C_m\cap X_W\bigr)$ lies in the first summand  $\pi_1(W) $.
Thus if $\pi_1(W) $ is infinite then even the  weak Lefschetz property
fails.

(\ref{3.key.exmps}.6) Continuing with the previous example, 
let $\tilde w:\tilde C_Q\to \tilde X$ be the family of conics
(that is rational curves of bidegree $(1,1)$) 
 and $w:C_Q\to X$ the 
corresponding family of conics in $X$.
Here both $\tilde w$ and $w$ have connected fibers.
It follows from our results that $\tilde w:\tilde C_Q\to \tilde X$
satisfies the  Lefschetz property but $w:C_Q\to X$
only satisfies the weak  Lefschetz property (with $N=2$).

(\ref{3.key.exmps}.7) The Lefschetz properties are really about
small open subsets of $X$ and of $C_M$. To see this,
let $u:Y\to X$ be a morphism between normal varieties,
 $X^0\subset X$ an open subvariety and $Y^0:=u^{-1}(X^0)$. Then
$\pi_1(X^0)\to \pi_1(X)$ is surjective
(cf.\ \cite[2.10]{shaf-book}), thus the index 
in Definition \ref{leff.prop.defn}
 increases as we pass from $X$ to $X^0$. That is,
$$
\Bigl| \pi_1(X^0) : 
\im\bigl[\pi_1(Y^0)\to \pi_1(X^0)\bigr]\Bigr|
\geq 
\Bigl|\pi_1(X) : \im\bigl[\pi_1(Y)\to \pi_1(X)\bigr]\Bigr|.
$$
Next let $C^0_M\subset C_M$ be a dense open subset. 
Then $C^0_m$ is a dense open subset of $C_m$ for general $m\in M$.
Thus, if $C_m$ is normal, then $\pi_1(C^0_m)\to \pi_1(C_m)$ is surjective
by \cite[2.10]{shaf-book}, so
$$
\im\bigl[\pi_1(C^0_m)\to \pi_1(X^0)\bigr]=
\im\bigl[\pi_1(C_m)\to \pi_1(X^0)\bigr].
$$
\end{exmp}

Our first result says that these examples  almost explain 
every failure of  the  Lefschetz property.

\begin{prop} \label{leff-ewkleff.prop}
Let  $X$ be  a   normal variety over $\c$ and
  $M\leftarrow C_M \to  X$   a
flat,  irreducible family of
 irreducible  varieties. Then each of the following
statements implies the next.
\begin{enumerate}
\item $M\leftarrow C_M\to X$
satisfies the  Lefschetz property. 
\item   $C_M\to X$ is dominant, has geometrically irreducible generic fiber
and $X$ is generically $C_M$-connected.
\item $M\leftarrow C_M\to X$
satisfies the  weak Lefschetz property. 
\end{enumerate}
\end{prop}

 In any concrete situation is  usually easy to check that 
$C_M\to X$ is dominant and has geometrically irreducible generic fiber.
Being generically $C_M$-connected is not always clear but it  holds if
$X$ is smooth, proper, has Picard number 1 and
$M\leftarrow C_M\stackrel{u}{\to} X$ is a 
locally complete family of free  curves; see \cite[IV.4.14]{rc-book}.

Sometimes the difference between the 
Lefschetz property and the  weak Lefschetz property
is minor, but in the arithmetic applications 
\cite{MR1786496, MR1745009, MR2019976}
having surjectivity is essential. The following main technical result
says that if we avoid the  bad situations (\ref{3.key.exmps}.1--3)
and we have surjectivity for $\pi_1(X)$ itself then
 the Lefschetz property holds. More generally, the extent of any failure of
the  Lefschetz property is determined by $X$ itself. 

\begin{thm} \label{main.thm}
Let $M\stackrel{p}{\leftarrow} C_M\stackrel{u}{\to} X$ be a family of
varieties over a  smooth (not necessarily proper) 
variety $X$, defined over $\c$.
 Assume that 
\begin{enumerate}
\item $p$ and $u$ are both smooth
with irreducible fibers,
\item $u$ is surjective and
\item $X$ is generically $C_M$-connected.
\end{enumerate}
Let $j:X^0\into X$ be an open subset  
and $j_*:\pi_1(X^0)\to \pi_1(X)$ the induced map on the fundamental groups.
 Then there is an open subset $\emptyset\neq M^0\subset M$ such that
$$
\im\bigl[\pi_1\bigl(C_m\cap u^{-1}(X^0)\bigr)\to \pi_1(X^0)\bigr]
=
j_*^{-1}\Bigl(\im\bigl[\pi_1\bigl(C_m\bigr)\to \pi_1(X)\bigr]\Bigr)
\eqno{(\ref{main.thm}.4)}
$$
for  every $m\in M^0$.
\end{thm}

We already know from Proposition \ref{leff-ewkleff.prop}
that both images in (\ref{main.thm}.4) are finite index subgroups.
Thus (\ref{main.thm}.4) is equivalent to the equality
$$
\Bigl| \pi_1(X^0) : 
\im\bigl[\pi_1\bigl(C_m\cap u^{-1}(X^0)\bigr)\to \pi_1(X^0)\bigr]\Bigr|
=
\Bigl|\pi_1(X) : \im\bigl[\pi_1(C_m)\to \pi_1(X)\bigr]\Bigr|.
$$
If $X$ is simply connected then the right hand side of
 (\ref{main.thm}.4)
equals $\pi_1(X^0)$. Thus, in this case, we assert that
$\pi_1\bigl(C_m\cap u^{-1}(X^0)\bigr)\to \pi_1(X^0)$ is onto
for every $m\in M^0$. The latter  is  exactly the Lefschetz property.

When applying Theorem \ref{main.thm} to any 
family $M\stackrel{p}{\leftarrow} C_M\stackrel{u}{\to} X$,
we first replace $M$ by $M\setminus \sing M$, then
replace $C_M$ by the largest open subset $C^0_M$ where 
$p$ and $u$ are both smooth
and finally replace $X$ by $u(C^0_M)$. The first step is
entirely harmless.  The key question is to understand
how large $X\setminus u(C^0_M)$ is; only the divisors
contained in it matter.
 
 As a significant example, let  $X$ be  smooth, proper 
and  $M\subset \mor(\p^1, X)$  a nonempty, irreducible, open subset
 with universal morphism 
 $u:M\times \p^1\to X$. For $x\in X$ let
$M_x\subset M$ be the set of maps $[f]\in M$ such that
$f(0{:}1)=x$.

\begin{cor} \label{RC.leff.cor}
Let $X$ be a normal, proper variety
and  $M\subset \mor(\p^1, X)$  a nonempty,
 irreducible, open subset
parametrizing free maps with universal morphism  $u:M\times \p^1\to X$.
Assume that
\begin{enumerate}
\item $X\setminus \sing X$ is simply connected,
\item  $X\setminus u\bigl(M\times \p^1\bigr)$ has codimension $\geq 2$ and
\item $X$ is generically $M\times \p^1$-connected.
\end{enumerate}
Then $u:M\times \p^1\to X$
satisfies the Lefschetz property iff $M_x$ is irreducible for
general $x\in X$.
\end{cor}

Proof. The projection $M\times \p^1\to M$ is obviously smooth and
$u$ is smooth by \cite[I.3.5.4]{rc-book} since we parametrize free morphism.
We apply Theorem \ref{main.thm} to
$X^*:=u\bigl(M\times \p^1\bigr)$ replacing $X$. By assumption, 
$X^*$ is obtained from the 
simply connected smooth variety $X\setminus \sing X$ by
removing a closed subscheme of codimension $\geq 2$.
Thus $X^*$ is also simply connected and hence
 the right hand side of (\ref{main.thm}.3)
equals $\pi_1(X^0)$. \qed
\medskip

\begin{rem} \label{bouquet.rem} If $M_x$ is reducible for
general $x\in X$ then instead of  $u:M\times \p^1\to X$
one can work with the family of rational curves obtained
by smoothing a bouquet of rational curves through $x$,
one from each irreducible  component of $M_x$.
\end{rem}

Note that the assumptions (\ref{RC.leff.cor}.1--3)  hold
if $X$ is smooth and has  Picard number
$\rho(X)=1$. Thus we get the following.

\begin{cor} \label{RC.leff.cor.cor}
Let $X$ be a smooth proper variety with $\rho(X)=1$. 
Let  $M\subset \mor(\p^1, X)$  be a nonempty,
 irreducible, open subset
parametrizing free maps. 
Then the  universal morphism $u:M\times \p^1\to X$
satisfies the Lefschetz property iff $M_x$ is irreducible for
general $x\in X$. \qed
\end{cor}

\begin{say}[Proof of Corollary \ref{lines.on.hyp.cor}]  
Let $M\leftarrow C_M\stackrel{u}{\to} X$
be the universal family of lines.
Let $x\in X$ be a point. After a coordinate change,
we may assume that $x=(1{:}0{:}\cdots{:}0)$.
Write the equation of $X_d$ as
$$
g_1(x_1,\dots, x_n)x_0^{d-1}+\cdots + g_d(x_1,\dots, x_n).
$$
The family of lines in $X$ through $x$ is then given by
the equations
$$
M_x:=\bigl(g_1=\cdots=g_d\bigr)\subset \p^{n-1}.
$$

$M_x$ is smooth  of dimension $n-1-d$ for general $x\in X$ by
\cite[II.3.11]{rc-book}. Thus $M_x$ is a smooth complete intersection, hence
irreducible if $n-1-d\geq 1$.

Thus $M$ has a unique irreducible component
$M^0\subset M$ such that the corresponding family
$u^0: C^0_M{\to} X$ is dominant and has
geometrically irreducible generic fiber.
Thus, by Corollary \ref{RC.leff.cor.cor},
 $M^0\leftarrow C^0_M\stackrel{u^0}{\to} X$
satisfies the Lefschetz property.

Conversely, assume that $d\geq n-1$. If 
 $d\geq n$ then there is no line through a general point;
this is like example (\ref{3.key.exmps}.1).
The $d=n-1$ case is discussed in (\ref{3.key.exmps}.4). \qed
\end{say}

\begin{rem} The  proof applies to any smooth, Fano complete
intersection of Fano index $\geq 3$. If the Fano index is $2$,
applying Remark \ref{bouquet.rem} yields very high degree curves,
but most likely
conics work if the Fano index is $2$ and cubics if 
 the Fano index is $1$.

As far as I know, Corollary \ref{lines.on.hyp.cor}
should hold in any characteristic. It holds for
general hypersurfaces where the family of lines is smooth and has
the expected dimension (cf.\ \cite[V.4.3]{rc-book}).
\end{rem}

The proof of Theorem \ref{main.thm} follows the outlines
of \cite[Sec.5]{k-hom}. First we recall properties of open chains,
then we pass to a subfamily that is topologically trivial.
After studying which chains lift to \'etale covers, the
proof is completed in Paragraphs \ref{all.lift.say}--\ref{key.lift.prop.2}.
At the end we consider how to modify the statements and the proofs
to work over arbitrary fields.

\subsection*{Open chains}{\ }

\begin{say}[Chains of varieties over $X$] \label{chains.say}
Let 
$M\leftarrow C_M\stackrel{u}{\to} X$  be a family of
schemes over $X$.
A  {\it $C_M$-link} is a morphism of a  triple $u_m:( C_m, a, b)\to X$
where $m\in M$ and $a,b\in C_m$. 
A {\it $C_M$-chain} of length $r$  over $X$ 
consists of
\begin{enumerate}
\item $C_M$-links  $u_i:(C_i, a_i, b_i)\to X$  for $i=1,\dots, r$ such that
\item $u_i(b_i)=u_{i+1}(a_{i+1})$
 for $i=1,\dots, r-1$.
\end{enumerate}
We say that the chain {\it starts} at $u_1(a_{1})\in X$ 
and {\it ends}
at $u_r(b_r)\in X$ or that it   {\it connects}
$u_1(a_{1})$ and $u_r(b_r)$.

A  $C_M$-chain determines a  reducible variety
$\vee_{i=1}^rC_i$
obtained from the  disjoint union of  $C_1,\dots, C_r$  by identifying
$b_i\in C_i$ with $a_{i+1}\in C_{i+1}$ for $i=1,\dots, r-1$.
The morphisms $u_i$ then define a morphism
$\vee_i u_i: \vee_iC_i\to X$.
If the $C_i$ are connected then the
 image of $\vee_i u_i $ is a connected subscheme of $X$ which contains
the starting and end points of the chain.

Starting with $M\leftarrow C_M\stackrel{u}{\to} X$
the set of all pairs  $( C_m,a)$ (resp.\ triples $( C_m,a, b)$)  
is naturally given by
$$
C_M\leftarrow C_M\times_MC_M\stackrel{u\circ \pi_2}{\longrightarrow} X
\qtq{and}
C_M\times_MC_M\leftarrow C_M\times_MC_M\times_MC_M\stackrel{u\circ \pi_3}
{\longrightarrow} X
$$
where the marked points are given by the diagonal maps
$$
\delta:C_M\to C_M\times_MC_M\qtq{and}
\delta_1, \delta_2:C_M\times_MC_M\to C_M\times_MC_M\times_MC_M.
$$
Here $\pi_i$ denotes the $i$th coordinate projection and
$\delta_i$ maps the first $C_M$ identically to the $i$th factor
and the second $C_M$ diagonally to the other 2 factors.
Thus all $C_M$-chains of length 1 are parametrized by 
$$
C_M\times_MC_M\leftarrow C_M\times_MC_M\times_MC_M\stackrel{u\circ \pi_3}
{\longrightarrow} X
\eqno{(\ref{chains.say}.3)}
$$
which we denote from now on by 
$$
\chain(C_M, 1)\stackrel{p^{(1)}}{\longleftarrow}
 C_M^{(1)}\stackrel{u^{(1)}}{\longrightarrow} X.
\eqno{(\ref{chains.say}.4)}
$$
Out of this we get that all $C_M$-chains of length 2 are
parametrized by 
$$
\chain(C_M, 2) :=\chain(C_M, 1)\times_X \chain(C_M, 1)
$$
where the 2 maps  $\chain(C_M, 1)\to X$ are given by
$u^{(1)}\circ \delta_2$ on the first copy and
$u^{(1)}\circ \delta_1$ on the second copy.
Over this there is a universal family
$$
\chain(C_M, 2)\ \stackrel{p^{(r,1)}\vee p^{(r,2)}}{\longleftarrow}\ 
 C_M^{(2,1)}\vee C_M^{(2,2)}\ \stackrel{u^{(r,1)}\vee u^{(r,2)}}{\longrightarrow}\  X
$$
where $ C_M^{(r,i)} $ denotes the universal family of the $i$th links of the
$r$-chains.

By iterating this we get $\chain(C_M, r)$ parametrizing
length $r$ chains
$$
\chain(C_M, r)\ \stackrel{\vee_i p^{(r,i)}}{\longleftarrow}\ 
 \vee_{i=1}^rC_M^{(r,i)}\ \stackrel{\vee_i u^{(r)}}{\longrightarrow}\  X.
\eqno{(\ref{chains.say}.5)}
$$

If $C_M\to M$ is flat with irreducible fibers then
$C_M^{(1)}\to \chain(C_M, 1)$ is also flat with irreducible fibers.
For a point $x\in X$ we have
$$
\bigl(u^{(1)}\bigr)^{-1}(x)\cong u^{-1}(x)\times_M C_M.
$$
Thus we conclude the following.
\medskip

{\it Claim \ref{chains.say}.6.} Assume that $M$ is irreducible and both
maps  $M\leftarrow C_M\to X$ are flat with irreducible fibers.  Then:
\begin{enumerate}
\item[a)]  Each $\chain(C_M, r)$ is irreducible.
\item[b)] The maps
$\chain(C_M, r)\stackrel{p^{(r,i)}}{\longleftarrow}
 C_M^{(r,i)}\stackrel{u^{(r,i)}}{\longrightarrow} X$ 
are flat with irreducible fibers.
\item[c)]   If $C^0_M\subset C_M$ is a dense open subset
then $\chain(C^0_M, r)$ is a dense open subset of $\chain(C_M, r)$.
\qed
\end{enumerate}
\end{say}

\begin{defn} \label{chains.say.8} With the above notation, 
the starting and end points
give  morphisms
$$
\alpha, \beta: \chain(C_M,r)\to X. 
$$
We say that $X$ is  {\it generically $C_M$-connected}
if 
$$
\alpha\times\beta: \chain(C_M,r)\to X\times X 
$$
is dominant for some $r$, that is, if two general points of $X$ can be 
connected by a $C_M$-chain of length $r$.
(The equivalence of this definition with the one given in (\ref{3.key.exmps}.2)
is proved in  \cite[IV.4.13]{rc-book}.)

If $u$ is open and $X$ is  generically $C_M$-connected
then $\alpha\times\beta $ is 
dominant for every $r\geq\dim X$ by \cite[IV.4.13]{rc-book}.

Note that if  $X$ is  generically $C_M$-connected, $M$ is irreducible and both
maps  $M\leftarrow C_M\to X$ are flat with irreducible fibers
 then, by  (\ref{chains.say}.6.c), 
$X$ is also  generically $ C_M^0$-connected
for every  dense open subset  $ C_M^0\subset C_M$.
\end{defn}

Now we choose an especially well behaved subset $ C_M^0\subset C_M$.

\begin{prop} \label{gen.pos.trans.prop} Let $X$ be a normal variety and
$M\stackrel{p}{\leftarrow} C_M\stackrel{u}{\to} X$
a family of varieties over $X$ 
where  $M$ is irreducible and both
maps   are flat with irreducible fibers.
Let $\emptyset\neq X^0\subset X$ be an open subset. 
Then there is an  open subset $\emptyset\neq C_M^0\subset C_M$
with induced maps $p^0:C_M^0\to M$ and $u^0:C_M^0\to X$
such that
\begin{enumerate}
\item $p^0$ is smooth with irreducible fibers,
\item $p^0$ is a  topologically locally trivial
fiber bundle (over its image),
\item the image of $u^0$ is contained in $X^0$ and
\item $u^0$ has irreducible (hence connected) fibers. 
\end{enumerate}
\end{prop}

Proof.
We first replace $C_M$ by the open subset $C_M^1=u^{-1}(X^0)$
and then  by the open subset $C_M^2\subset C_M^1$ where $p$ is smooth. 

By  \cite[p.43]{gm-book}, 
every map between algebraic varieties is a 
locally topologically trivial fiber bundle
over a Zariski open subset. Thus by passing to an 
open subset
$C_M^0\subset C_M^2$ we may assume that 
properties (1--3) hold. Since each fiber of $u$ is irreducible,
the same holds for $u^0$.
\qed
\medskip

The pointed fibers  $(C^0_m, a)$ also form a
topologically locally trivial
fiber bundle  $C^0_M\leftarrow C_M^0\times_MC_M^0$.
Given a point $x\in X^0$, the set of all   $(C^0_m, a)$
such that $u(a)=x$ form a topologically locally trivial
fiber bundle over the connected base $(u^0)^{-1}(x)$.
As we noted in  (\ref{3.key.exmps}.7),
$$
\im\bigl[\pi_1(C^0_m,a)\to \pi_1(X^0,x)\bigr]=
\im\bigl[\pi_1(C^1_m,a)\to \pi_1(X^0,x)\bigr].
$$
Thus we obtain the following.

\begin{cor}\label{pi1.triv.subfamily} Notation and assumptions as in
(\ref{gen.pos.trans.prop}).
Then, for every $m\in M$, $a\in C^0_m$ and $x:=u(a)$, 
the image of the induced map
 $$
\Gamma(X^0,C,x)
:=\im\bigl[\pi_1\bigl(C^0_m, a\bigr)
\stackrel{u^0_*}{\longrightarrow} \pi_1\bigl( X^0, x\bigr)\bigr]
\subset \pi_1\bigl( X^0, x\bigr)
 $$
depends only on $\bigl( X^0, x\bigr)$ and $C_M$ but 
 not on $m\in M$ and $a\in C^0_m$.\qed

\end{cor}

\subsection*{Topologically locally trivial  families}

\begin{say} \label{gen.pos.trans.say}
We work  with 
families $M\stackrel{p}{\leftarrow} C_M\stackrel{u}{\to} X$ such that
$p$ has irreducible fibers and the following holds:

(\ref{gen.pos.trans.say}.1) For every $x\in X, m\in M$ and $a\in C_m$ satisfying
$u_m(a)=x$,  the image of the induced map
 $$
\Gamma(X,C,x)
:=\im\bigl[\pi_1\bigl(C_m, a\bigr)
\stackrel{u_*}{\longrightarrow} \pi_1\bigl( X, x\bigr)\bigr]
\subset \pi_1\bigl( X, x\bigr)
 $$
does not depend on  $m\in M$ and $a\in C_m$.

An equivalent formulation is the following. 

(\ref{gen.pos.trans.say}.2) Let $\bigl(\tilde X, \tilde x\bigr)\to (X,x)$
be any covering space such that 
$$
u_m: \bigl(C_m,a\bigr)\to \bigl(X,x\bigr)
\qtq{lifts to} 
\tilde u_m: \bigl(C_m,a\bigr)\to 
\bigl(\tilde X,\tilde x\bigr)
$$
for some  $m\in M$ and $a\in C_m$. Then
 the lift exists
for every $n\in M, b\in C_n$  for which
$u_n(b)=x$.

Now fix a point $x\in X$. Corresponding to
$\Gamma(X,C,x)$
there is an \'etale cover
$$
q_X: \bigl(\tilde X, \tilde x\bigr)\to \bigl(X,x\bigr).
\eqno{(\ref{gen.pos.trans.say}.3)}
$$
We do not yet know that $\Gamma(X,C,x)$ has finite index,
so $\tilde X\to X $ could have infinite degree.
Thus $\tilde X$ is an analytic space for now.
\end{say}

\begin{prop}\label{key.lift.prop.1}
Notation and assumptions as in (\ref{gen.pos.trans.say}). Then 
every $C_M$-chain on $X$ starting at $x$ lifts to a 
$C_M$-chain on $\tilde X $  starting at $\tilde x$.
\end{prop}

Proof. A  $C_M$-chain is given by the data
$u_i:(C_i, a_i, b_i)\to X$. Set $x_1:=x$. By the choice of $\Gamma(X,C,x_1)$,
$$
 u_1: (C_1, a_1)\to (X, x_1)
\qtq{lifts to}
\tilde u_1: (C_1, a_1)\to (\tilde  X, \tilde x_1).
$$
If we let $\tilde x_2$ denote the image of $b_1$ 
then we can view the latter map as
$$
\tilde u_1: (C_1, b_1)\to (\tilde  X, \tilde x_2).
$$
We next apply (\ref{gen.pos.trans.say}) to 
$$
u_1: ( C_1, b_1)\to (X, x_2)\qtq{and}
u_2: ( C_2, a_2)\to (X, x_2)
$$
to see that if one of them lifts to 
$(\tilde  X, \tilde x_2)$ then so does the other.
This gives us
$$
\tilde u_2: (C_2, a_2)\to (\tilde  X, \tilde x_2).
$$
We can iterate the argument to lift the whole chain. \qed
\medskip

\begin{cor}\label{key.lift.cor}
Notation and assumptions as in (\ref{gen.pos.trans.say}).
Assume in addition that $X$ is $C_M$-connected. Then
$\Gamma(X,C,x)\subset  \pi_1\bigl(X,x\bigr)$, as in (\ref{gen.pos.trans.say}.1),
has finite index, thus 
$q_X: \bigl(\tilde X, \tilde x\bigr)\to \bigl(X,x\bigr)$
is an algebraic  \'etale cover.

More precisely, the degree of  $q_X$ is bounded by $N:=$ the
number of irreducible components of the geometric generic fiber of
$\alpha\times\beta:\chain(X,\dim X)\to X\times X$.
\end{cor}

Proof. Let $\chain(C_M,r,x)\subset \chain(C_M,r)$ denote the
subscheme parametrizing chains that start at $x$. Thus
$\chain(C_M,r,x)$ is a fiber of 
$\alpha: \chain(C_M,r)\to X$ and, for general $x\in X$,
the number of irreducible components of the geometric generic fiber of
$\beta:\chain(X,r,x)\to X$ equals $N$.

Let $p^{(r)}:\vee_i C^{(r,i)}_M\to \chain(X,r,x)$ be the universal family
of $C_M$-chains of length $r$ with starting and end point sections
$\alpha, \beta: \chain(X,r,x)\to \vee_i C^{(r,i)}_M$. 

Note that  $u^{(r,1)}\circ \alpha$ maps $\chain(X,r,x)$ to $\{x\}$
and,
by (\ref{key.lift.prop.1}),
each fiber of  $p^{(r)}$ lifts to a $\tilde C_m$-chain on $\tilde X$
starting at $\tilde x$. Thus
$$
\vee_i u^{(r,i)}: \vee_i C^{(r,i)}_M\to X
\qtq{lifts to}
\vee_i \tilde u^{(r,i)}: \vee_i \tilde C^{(r,i)}_M\to \tilde X.
$$
 In particular, the end point map
$$
u^{(r,r)}\circ \beta: \chain(C_M,r,x)\to X
\qtq{lifts to} 
\tilde u^{(r,r)}\circ \tilde \beta: \chain(C_M,r,x)\to \tilde X.
$$
Therefore
$$
\im\bigl[ \beta_*:\pi_1\bigl(\chain(C_M,r,x)\bigr)\to \pi_1(X)\bigr]
\subset \Gamma(X,C,x).
$$
By assumption (and \cite[4.13]{rc-book}) $\beta$ is dominant
for $r\geq \dim X$. Therefore, by \cite[2.10]{shaf-book},
  the index is bounded as 
$$
\Bigl|\pi_1(X) : \im\bigl[ \pi_1\bigl(\chain(C_M,r,x)\bigr)\to \pi_1(X)
\bigr]\Bigr|\leq N. \qed
$$

\begin{say}[Proof of \ref{leff-ewkleff.prop}]
The implication (\ref{leff-ewkleff.prop}.1) $\Rightarrow$
(\ref{leff-ewkleff.prop}.2)
was already noted in (\ref{3.key.exmps}.1--3).

It remains to show that (\ref{leff-ewkleff.prop}.2) $\Rightarrow$
(\ref{leff-ewkleff.prop}.3).

As we noted in (\ref{chains.say.8}),
replacing $C_M$ with an open subset $\emptyset\neq C^0_M\subset C_M$
does not change the assumptions. Thus we may assume that
the assumptions of (\ref{key.lift.cor}) hold.
This gives the bound $N:=$ the
number of irreducible components of the geometric generic fiber of
$\alpha\times\beta:\chain(X,\dim X)\to X\times X$. \qed

\end{say}

\subsection*{Proof of Theorem \ref{main.thm}}{\ }

Fix an open subset $X^0\subset X$ and use
(\ref{gen.pos.trans.prop}) to obtain $C^0_M\subset C_M$. Then pick a general
point $x\in X^0$ and, as in Paragraph \ref{gen.pos.trans.say},
construct
$$
q^0_X: \bigl(\tilde X^0, \tilde x\bigr)\to (X^0,x).
$$
By Proposition \ref{leff-ewkleff.prop},   $q^0_X $
has finite degree, thus it extends (uniquely)
to a normal, possibly ramified, finite cover
$$
q_X: \bigl(\tilde X, \tilde x\bigr)\to (X,x).
$$
If $q_X$ is also \'etale then $\tilde X^0\to X^0 $
is the pull-back of the finite \'etale cover
$\tilde X\to X$; this is what (\ref{main.thm}.4) asserts.

All that remains is to derive a contradiction
if $q_X$ is ramified.
Since $X$ is smooth,  in this case there is 
 a nonempty branch divisor  $B\subset X$. 
First we show that most $C_M$-chains starting at $x$ lift to
 $\tilde X$. Then we use the branch divisor to show that most
chains do not lift, thereby arriving at a contradiction.

\begin{say}[Lifting  $C_M$-chains]\label{all.lift.say}
 A  $C^0_M$-chain is given by the data
$u^0_i:(C^0_i, a_i, b_i)\to X$. Here each $C^0_i $ is an open subset
of the corresponding $C_i$ thus  the $C^0_M$-chain
naturally corresponds to a  $C_M$-chain given by the data
$u_i:(C_i, a_i, b_i)\to X$.

 Since the $C_i$ are normal (even smooth)
and $q_X:\tilde X\to X$ is finite, every lifting 
$\tilde u^0_i: C^0_i\to \tilde X$ of $u^0_i$ uniquely extends
to $\tilde u_i: C_i\to \tilde X$.  Thus if a $C^0_M$-chain
lifts to $\tilde X^0 $ then the corresponding  $C_M$-chain also
lifts to $\tilde X$.
\end{say}

\begin{say}[Non-liftable chains] \label{key.lift.prop.2}
Assume that the branch divisor $B_X\subset X$ of 
 $q_X:\tilde X\to X$ is nonempty.
Let $B^*_X\subset B_X$ be the open subset of smooth points. 
Since $u:C_M\to X$ is surjective and smooth, the preimage
$B_C:=u^{-1}(B_X)$ is also  nonempty and 
$u^{-1}(B^*_X)$ is smooth and nonempty. 
Let $B^*_C\subset u^{-1}(B^*_X)$ be the set of points where
the restriction of $p$ to $B_C$ is smooth.
Finally let $M^*\subset M$ be the open subset consisting of
those points $m\in M$ such that $C_m$ meets $B^*_C$
in at least 1 point. Thus for $m\in M^*$ there is a map of the
unit  disc
$\tau_m:\Delta\to C_m$ such that
$u\circ \tau_m: \Delta\to X$ is transversal to $B$. 
Since $q_X$ branches along $B$, the sheets of $\tilde X\to X$
have nontrivial monodromy around $B$ and the pull-back
to $\Delta$ still has nontrivial monodromy.

 Set $d:=\deg \tilde X/ X$.
If $m\in M^*$ then the pull-back
$$
q_m: C_m\times_X\tilde X\to C_m
$$
is a degree $d$  cover that is \'etale outside  $C_m\cap B^*_C$ 
and whose  monodromy around $C_m\cap B^*_C$  is nontrivial. 
 The cover need not be connected or normal, but,
due to the monodromy, 
it can not be a  union of $d$ trivial covers
$C_m\cong C_m$.
That is, if  $a\in C_m$ is a general  point
and $\tilde a_1,\dots, \tilde a_d$ its preimages
in $C_m\times_X\tilde X $ then, 
 for at least one $\tilde a_i$, the identity map
$(C_m, a)\to (C_m, a)$ can not be lifted to
$( C_m, a)\to\bigl( C_m\times_X\tilde X, \tilde a_i\bigr)$.

Thus if $y\in X$ is the image of $a$ and
$\tilde y_1,\dots, \tilde y_d\in \tilde X$ its preimages,
then for at least one $\tilde y_i$, the map
$ u_m: (C_m, a)\to (X,y)$ can not be lifted to
$$
\tilde u_{(m,i)}:(C_m , a)\not\to 
\bigl( \tilde X, \tilde y_i\bigr).
$$
Consider now the dominant map
$\tilde \beta_r: \chain(C^0_M,r,x)\to \tilde X^0$
and let $X^*\subset X$ be a Zariski open subset
such that $q_X^{-1}(X^*)\subset \im \tilde \beta_r$.

By choosing the above $u_m:C_m\to X$
generically, we may assume that there is a  point
$a_{r+1}\in C_m$ such that $y:=u_m(a_{r+1})\in X^*$.

By the choice of $X^*$, for every $\tilde y^*_i\in q_X^{-1}(y)$
there is a $C^0_M$-chain of length $r$
whose lift to $\tilde X$ connects $\tilde x$ and $\tilde y^*_i$. 
We can add $ u^0_m: \bigl(C^0_m, a_{r+1}, b_{r+1}\bigr)\to X$ 
as the last link of
any of these chains.  Thus we get $d$ different 
$C^0_M$-chains of length $r+1$
and, for  at least one of them,
its extension to a $C_M$-chain
 can not be lifted to $\tilde X$.
This contradicts (\ref{all.lift.say}) and completes the proof of
 Theorem \ref{main.thm}. 
\qed

\end{say}

\subsection*{Other fields}{\ }

Our results apply to varieties over an arbitrary field,
with two modifications.

First, we have to use the algebraic fundamental group;
which we still denote by $\pi_1$.
Note that  if $k$ is any field with algebraic closure $\bar k$ and
$p:Y\to X$ is a morphism of geometrically irreducible $k$-varieties
then  the induced map $\pi_1(Y)\to \pi_1(X)$ is surjective iff
$\pi_1(Y\times_k\bar k)\to \pi_1(X\times_k\bar k)$ is surjective.
Thus our questions are geometric in nature and the 
key point is to understand what happens over algebraically closed
fields in positive characteristic.

The main difference is that even the classical Lefschetz theorem fails
in the non-projective case. For instance,
$\pi_1(\a^1)\to \pi_1(\a^2)$ is not surjective in positive characteristic.
(An example is given by the cover  $(z^p+z+x=0)\subset \a^3$ of
the $xy$-plane which splits over any line $x=c$.) 
This is remedied with the following variant of Definition \ref{leff.prop.defn}.

\begin{defn}\label{leff.prop.defn.gen}
 Let $k$ be an algebraically closed field  of positive characteristic
and
$$
M\stackrel{p}{\leftarrow} C_M\stackrel{u}{\to} X.
\eqno{(\ref{leff.prop.defn.gen}.1)}
$$
a family of schemes where  $M$ is geometrically irreducible
and $p$ is flat with irreducible  fibers.
We say that  the  family   (\ref{leff.prop.defn.gen}.1)
satisfies the {\it Lefschetz property} if the following holds.
\begin{enumerate}
\item[]  For every Zariski open dense subset $X^0\subset X$
and every finite quotient $\pi_1(X^0)\onto G$ 
there is a Zariski open dense subset $M^0_G\subset M$ such that, 
for every  $m\in M^0_G$, the induced map
$$
u(X^0,G,m)_*:\pi_1\bigl(C_m\cap u^{-1}(X^0)\bigr)\to 
\pi_1(X^0)\to G\qtq{is surjective.}
$$
\end{enumerate}
We say that  (\ref{leff.prop.defn.gen}.1)
satisfies the {\it weak Lefschetz property} if there is a constant
$N$  (independent of $X^0$  and of $G$) such that,
for a suitable choice of $ M^0_G$,
the image of  $u(X^0,G,m)_* $
has index at most $N$ in  in $G$.

\end{defn}

With this notion, the only question is what should
replace the topologically trivial family used in (\ref{gen.pos.trans.prop}). 
Topological triviality is used only through
its consequence    (\ref{gen.pos.trans.say}.1).
In our case we need that
 $$
\Gamma(X,C,x)
:=\im\bigl[\pi_1\bigl(C_m, a\bigr)\to \pi_1\bigl( X, x\bigr)\to G\bigr]
\subset G
 $$
be independent of $m\in M$ and $a\in C_m$.
 This is an easy consequence of the
semicontinuity property of the fundamental groups in fibers;
see  \cite[Prop.16]{MR2011744} for a precise statement and  proof.

The rest of the arguments go through with minor changes.

 \begin{ack}
I thank J.~Starr
 for useful comments. 
Partial financial support   was provided  by  the NSF under grant number 
DMS-0968337.
\end{ack}

\def\cprime{$'$} \def\cprime{$'$} \def\cprime{$'$} \def\cprime{$'$}
  \def\cprime{$'$} \def\cprime{$'$} \def\dbar{\leavevmode\hbox to
  0pt{\hskip.2ex \accent"16\hss}d} \def\cprime{$'$} \def\cprime{$'$}
  \def\polhk#1{\setbox0=\hbox{#1}{\ooalign{\hidewidth
  \lower1.5ex\hbox{`}\hidewidth\crcr\unhbox0}}} \def\cprime{$'$}
  \def\cprime{$'$} \def\cprime{$'$} \def\cprime{$'$}
  \def\polhk#1{\setbox0=\hbox{#1}{\ooalign{\hidewidth
  \lower1.5ex\hbox{`}\hidewidth\crcr\unhbox0}}} \def\cdprime{$''$}
  \def\cprime{$'$} \def\cprime{$'$} \def\cprime{$'$} \def\cprime{$'$}
\providecommand{\bysame}{\leavevmode\hbox to3em{\hrulefill}\thinspace}
\providecommand{\MR}{\relax\ifhmode\unskip\space\fi MR }
\providecommand{\MRhref}[2]{%
  \href{http://www.ams.org/mathscinet-getitem?mr=#1}{#2}
}
\providecommand{\href}[2]{#2}

\vskip1cm

\noindent Princeton University, Princeton NJ 08544-1000

{\begin{verbatim}kollar@math.princeton.edu\end{verbatim}}

\end{document}